\documentclass{article}
\usepackage[all]{xy}
\usepackage{amsfonts,amsmath,oldgerm,amssymb,amscd,theorem}

\newcommand{\surj}{\twoheadrightarrow}
\newcommand{\inj}{\hookrightarrow}
\newcommand{\lra}{\longrightarrow}
\newcommand{\Q}{\mbox{$\mathbb Q$}}     
\UseComputerModernTips

\newcommand{\ov}{\overline}

\newtheorem{theorem}{Theorem}[section]

\newtheorem{lemma}[theorem]{Lemma}

\newtheorem{corollary}[theorem]{Corollary}

\newtheorem{sublemma}[theorem]{Sublemma}
\newenvironment{conjecture}{\newline{\bf Conjecture.}\it}{\newline}

\theoremstyle{plain}
\theorembodyfont{\upshape}

\begin{document}
\title{On $K_2$ of 1-dimensional local rings}
\author{Amalendu Krishna}
\date{}
\baselineskip=14pt   
\maketitle                                                 
\begin{abstract}
We study $K_2$ of one-dimensional local domains over a field of
characterstic 0, introduce a conjecture, and show that this conjecture
implies Geller's conjecture. We also show that Berger's conjecture
implies Geller's conjecture, and hence verify it in many new cases. 
\end{abstract}
\section{Introduction}
Let $A$ be an one-dimensional local integral domain which is essentially
of finite type over a field $K$ of characteristic zero. Let $F$ denote
the field of fractions of $A$. It is easy to check that the map $K_1(A)
\lra K_1(F)$ is injective. It was a general question if the map
$K_2(A) \lra K_2(F)$ is also injective, which is now known not to be true
in general. Quillen's proof of Gersten's conjecture implies that this
is true if $A$ is a regular local ring. Dennis and Sherman ([G]) showed
that this map is not injective when $A$ is the local ring of the singular
point of the cuspidal curve Spec$(K[t^2, t^3])$. The general picture
about singular rings is given by the following conjecture of Geller ([G]).
\\
\begin{conjecture} ({\bf Geller})
Let $A$ be a local one-dimensional domain with field of fractions $F$. Then
$A$ is regular if and only if the map $K_2(A) \lra K_2(F)$ is injective.
\end{conjecture}

This conjecture was verified by Geller ([G]) when $A$ is noetherian, 
equicharacteristic, characteristic zero, and is also seminormal with finite
normalisation. In the same article, Dennis and Sherman verify it for cuspidal
rings of the type $K[t^2, t^3]$ as described above. The conjecture is
still unknown in almost all other cases. Our first aim in this paper is to
formulate an Artinian version of this conjecture, and to show that this
conjecture implies Geller's conjecture. Before we state the conjecture,
let us recall that that an algebra over a field $K$ is called a
$principal \ ideal \ algebra$ if every ideal of $A$ is principal. We
call $A$ to be finite-dimensional if it is finite over $K$. In this
paper, we will have standing assumption that $K$ is an algebraically 
closed field of characteristic zero, and all $K$-algebras are 
essentially of finite type over $K$. Our Artinian version of above 
conjecture is
\begin{conjecture} ({\bf AGC}) If $A$ is a subalgebra of a finite-dimensional
principal ideal $K$-algebra $B$ such that the map $K_2(A) \lra K_2(B)$
is injective, then $A$ is also a principal ideal algebra.
\end{conjecture}
\\
\\      
We shall call this `Artinian Geller Conjecture' (AGC).
Our first result in this paper states
\begin{theorem}\label{thm1}
With $A$ and $K$ as above, Artinian Geller Conjecture implies Geller's 
Conjecture.
\end{theorem}
In the other part of this paper, our aim is to relate these conjectures
with differential forms in order to verify Geller's conjecture in
some cases. In this regard, we recall a similar conjecture about the
module of Kahler differentials on one-dimensional local domain over
a field of char. 0. 
\begin{conjecture} ({\bf Berger}) Let $A$ be an one-dimensional local domain
which is essentially of finite type over a field $K$ of characteristic zero. 
Then $A$ is regular if and only if the module of Kahler differentials 
${\Omega}_{A/K}$ is torsion-free.
\end{conjecture}
\\
\\
This conjecture was formulated by R. Berger in [B1] almost forty years ago.
This has been verified in many cases (listed below) by various people 
though it is still unknown in general.
\begin{theorem}\label{thm2}
If $A$ is as above, with $K$ algebraically closed, then Berger's Conjecture
implies Geller's Conjecture.
\end{theorem}

Before we state our corollary to this theorem, we recall that ([B2]) a local
ring $A$ as above is called an `almost' complete intersection if the
first quadratic transform of $A$ is a complete intersection. The common
examples are local rings of plane curves or a curve through a smooth
point of a surface.    
\begin{corollary}\label{cor}
Geller's Conjecture is true in each of the following cases. \\
(i) $A$ is seminormal (also proved by Geller), \\
(ii) ${{\mathfrak M}^3}B \subset A$, where $B$ is the normalisation
of $A$ with Jacobson radical $\mathfrak M$, \\ 
(iii) $A$ is a complete intersection, \\
(iv) $A$ is almost complete intersection, \\
(v) $A$ is the local ring of the vertex of an 1-dimensional graded ring 
with vertex as only singular point, \\
(vi) $A$ has analyitically smoothable curve singularities, \\
(vii) $A$ has multiplicity $< \binom{m}{2}$, where $m$ is the embedding 
dimension of $A$, and \\
(viii) $A$ has deviation $\le 3$.
\end{corollary}
{\bf Remark.} We mention here that the condition of the field $K$ being
algebraically closed is only a technical one and one can reduce the 
general case to this case using the techniques of [G] and [CGW]. In fact,
it is shown in [CGW] that one can always assume $K$ to be algebraically
closed to prove Berger's conjecture.
\section{Some results on Hochschild and Cyclic homology}
In this section, we aim to prove some results concerning
Hochschild and cyclic homology of rings. We refer the reader to [LO]
for basic notions of Hochschild and Cyclic homology of rings. Let
$k$ be field of characteristic 0, and we assume all $k$-algebras
to be commutative. For any $k$-algebra maps $A \lra B$, Loday ([LO])
also defines the relative Hochschild homology $HH^k_*(A, B)$ over $k$
as the homology groups of the chain complex $Cone(C_{\bullet}(A) \lra
C_{\bullet}(B))$, where $C_{\bullet}(A)$ denotes the Hochschild complex
of $A$ etc. For an ideal $I$ of $A$, $HH^k_*(A, I)$ will be the relative
homology of $A$ and $A/I$. One defines relative Cyclic homology in 
similar way by taking the cone over the total cyclic complexes of the
two algebras. We also have the notion of relative $K$-theory as defined,
for example in [CS]. There are Chern class maps $K_i(A) \lra HH^k_i(A)$,
(Dennis trace maps) and by functoriality of fibrations of $K$-theory 
spectra and Hochschild
homology, one also has Chern class maps from relative $K$-theory to relative
Hochschild homology ([LO]), which are compatible with long exact
sequence of relative $K$-theory and Hochschild homology. It is known
that there are natural maps ${\Omega}^i_{A/k} \lra HH^k_i(A)$ and
$HH^k_i(A) \lra {\Omega}^i_{A/k}$ such that the composite is 
multiplication by $i!$. In particular, $HH^k_1(A)$ is same as the module
of Kahler differentials on $A$ over $k$. For an ideal $I$ of $A$,
let ${\Omega}^1_{(A, I)/k} : = {\rm Ker}({\Omega}^1_{A/k} \surj 
{\Omega}^1_{{A/I}/k})$. We begin with the following
\begin{lemma}\label{HH1}
Let $A$ be a $k$-algebra which is reduced and and is essentially of finite
type over $k$. Let $B$ be the normalisation of $A$, and let $I$
be a conducting ideal for this normalisation. Then, for all sufficiently large
$n$, the map
$$HH^k_1(A, I^n) \lra HH^k_1(B, I^n)$$ is injective.
\end{lemma}      
{\bf Proof.}
We use the following commutative diagram of exact sequences.
$$
\xymatrix@C.6pc{
0 \ar[r] & {\frac{HH^k_2(A/{I^n})} {HH^k_2(A)}} \ar[r] \ar[d] &
HH^k_1(A, I^n) \ar[r] \ar[d] & {\Omega}^1_{(A, I)/k} \ar[r] \ar[d] & 0 \\
0 \ar[r] & {\frac{HH^k_2(B/{I^n})} {HH^k_2(B)}} \ar[r] &
HH^k_1(B, I^n) \ar[r] & {\Omega}^1_{(B, I)/k} \ar[r] & 0 
}
$$
$${\rm Diagram} \ 1$$
It is enough to show that the vertical maps on the ends are injective.  
Put ${\ov {{\Omega}^1_{A/k}}} : = {\rm Ker}({\Omega}^1_{A/k} \lra
{\Omega}^1_{B/k})$. This module is supported on $V(A/I)$ and hence is
annihilated by $I^n$ for $n >> 0$. Thus,
\begin{equation}\label{eqn1} 
I^n({\ov {{\Omega}^1_{A/k}}}) = 0 \ {\rm for} \ n >> 0.
\end{equation}
Furthermore, since ${\Omega}^1_{(A, I^n)/k} = I^n{\Omega}^1_{A/k} +
d(I^n)$, one has a diagram of exact sequences
$$
\xymatrix@C.6pc{
0 \ar[r] & d(I^n) \ar[r] \ar@{=}[d] &  {\Omega}^1_{(A, I^n)/k} \ar[r] \ar[d] &
{\frac{I^n{\Omega}^1_{A/k}} {d(I^n) \cap {I^n{\Omega}^1_{A/k}}}} \ar[r]
\ar[d] & 0 \\
0 \ar[r] & d(I^n) \ar[r] &  {\Omega}^1_{(B, I^n)/k} \ar[r] &
{\frac{I^n{\Omega}^1_{B/k}} {d(I^n) \cap {I^n{\Omega}^1_{B/k}}}} \ar[r]
& 0.
}
$$  
This gives
$${\rm Ker}({\Omega}^1_{(A, I^n)/k} \lra  {\Omega}^1_{(B, I^n)/k})
 = {\rm Ker} ({\frac{I^n{\Omega}^1_{A/k}} {d(I^n) \cap 
{I^n{\Omega}^1_{A/k}}}} \lra
{\frac{I^n{\Omega}^1_{B/k}} {d(I^n) \cap 
{I^n{\Omega}^1_{B/k}}}}).$$
Next, we \\
{\bf Claim.} For $n >> 0$, $I^n{{\Omega}^1_{A/k}} \inj 
I^n{{\Omega}^1_{B/k}}$. \\
To prove the claim, notice that ${\Omega}^1_{A/k}$ is a finitely generated
$A$-module, and hence by Artin-Rees theorem, there exists $c > 0$ such that
for all $n > c$,
$$(I^n{{\Omega}^1_{A/k}} \cap {\ov {{\Omega}^1_{A/k}}}) \subset
I^{n-c}(I^c{{\Omega}^1_{A/k}} \cap {\ov {{\Omega}^1_{A/k}}}).$$
In particular, we get
$$(I^n{{\Omega}^1_{A/k}} \cap {\ov {{\Omega}^1_{A/k}}}) \subset
I^{n-c}({\ov {{\Omega}^1_{A/k}}}) = 0 \ {\rm for} \ n >> 0 \
{\rm by} \ {\rm {~\ref{eqn1}}},$$
which proves the claim. \\
Using this claim, we obtain
$${\rm Ker} ({\frac{I^n{\Omega}^1_{A/k}} {d(I^n) \cap 
{I^n{\Omega}^1_{A/k}}}} \lra
{\frac{I^n{\Omega}^1_{B/k}} {d(I^n) \cap 
{I^n{\Omega}^1_{B/k}}}}) = 
{\frac{I^n{\Omega}^1_{A/k} \cap (d(I^n) \cap 
{I^n{\Omega}^1_{B/k}})} {d(I^n) \cap {I^n{\Omega}^1_{A/k}}}}$$
$$ \ \ \ \ \ \ \ \ \ \ \ \ \ \ \  
= {\frac {d(I^n) \cap {I^n{\Omega}^1_{A/k}}} 
{d(I^n) \cap {I^n{\Omega}^1_{A/k}}}} = 0.$$
This, together with ~\ref{eqn1} implies that the right-most vertical map
in Diagram~1 is injective for all sufficiently large $n$.

Before we start proving the injectivity of the vertical map on the left, 
we make the convention that all Hochschild and cyclic homologies will
be considered over the given base field $k$ in the remaining part of this
lemma, and we will suppress this field $k$. We use the 
Hodge decomposition (or $\lambda$-decomposition)
([LO] or [C]) on Hochschild homology to get
$${\frac {HH_2(A/{I^n})} {HH_2(A)}} =
{\frac {HH^{(1)}_2(A/{I^n})} {HH^{(1)}_2(A)}} \oplus
{\frac {HH^{(2)}_2(A/{I^n})} {HH^{(2)}_2(A)}}.$$
But, for any $k$-algebra $A$, one has
$HH^{(1)}_2(A) = D^k_1(A)$, and $HH^{(2)}_2(A) = {\Omega}^2_{A/k}$ by [LO]
(chapter 4), where $D^k_*(A)$ denotes the Andre-Quillen homology of $A$
over $k$. Moreover, for any ideal $I \subset A$, the map
${\Omega}^2_{A/k} \lra {\Omega}^2_{{(A/I)}/k}$ is surjective. Thus, we
get ${\frac {HH_2(A/{I^n})} {HH_2(A)}} = 
{\frac {D^k_1(A/{I^n})} {D^k_1(A)}}$ and similarly for $B$.
Now, from [LO] (chapter 3), we have a diagram of exact sequences   
$$
\xymatrix@C.6pc{
D^k_1(A) \ar[r] \ar[d] & D^k_1(A/{I^n}) \ar[d] \ar[r] &
D^A_1(A/{I^n}) \ar[d] \\
D^k_1(B) \ar[r] & D^k_1(B/{I^n}) \ar[r] &
D^B_1(B/{I^n}),
}
$$
and $D^A_1(A/{I^n}) = I^n/{I^{2n}} = D^B_1(B/{I^n})$. This proves the required
injectivity.
$\hfill \square$
\\  
\begin{lemma}\label{HH2}
Let $A$ be a regular ring which is essentially of finite type over a field
$K$ of char. zero. Let $I \subset A$ be an invertible ideal. Then, for
any subfield $k \subset K$, and any $n \ge 0$, the natural map
$${\frac {HH^k_n(A/{I^2})} {HH^k_n(A)}} \lra 
{\frac {HH^k_n(A/{I})} {HH^k_n(A)}}$$
is zero.
\end{lemma}
{\bf Proof.}
Since Hochschild homology commutes with localisation, we can assume that
$R$ is a regular local ring and $I = (t)$ is a principal ideal. Let
$A \mapsto D^k_*(A)$ denote the Andre-Quillen homology functor. Then, these
homology groups are given by
$$D^k_*(A) : = H_*({\mathbb L}_{A/k}),$$ 
where ${\mathbb L}_{A/k}$ denotes the cotangent complex of $A$ over $k$
([LO]). We first claim that $D^k_i(A) = 0$ for $i > 0$, and
$D^k_i(A/I) = 0$ for $i > 1$. \\
First, notice that since $k$ is of char. zero, $D^k_i(K)$ is the direct
limit of $D^k_i(L)$, where $L$ is a subfield of $K$ and finitely generated
over $k$. Moreover, $L$ can be viewed as a finite extension of a purely
transcendental extension (of finite degree) over $k$. But, the Andre-
Quillen homology of finite extension vanishes in char. 0, and a purely
transcendental extension of finite degree is a localisation of a 
polynomial ring over $k$ for which the Andre-Quillen Homology again vanishes.
Since Andre-Quillen homology commutes with direct limits ([Q]),
we conclude that $D^k_i(K) = 0$ for $i > 0$.
Now, we use this fact and the exact sequence ([LO])
$$D^k_i(K) \otimes A \rightarrow D^k_i(A) \rightarrow D^K_i(A) \rightarrow
D^k_{i-1}(K) \otimes A,$$ to see that it is enough to prove the claimed
statement over $K$.
However, since $A$ is smooth over $K$, and $I$ is a local complete 
intersection ideal in $A$, we have 
$D^K_i(A) = 0$ for $i > 0$ and $D^K_i(A/I) = 0$ for $i > 1$ by the results
of Avramov and Halperin ([AH]). 

Since $A$ is smooth over $K$, there is an isomorphism $HH^k_i(A) \cong
{\Omega}^i_{A/k}$ for any subfield $k \subset K$ ([C]). Furthermore,
since $D^k_i(A/I) = 0$ for $i > 1$, we have by [LR] (Theorem~3.1 and
Proposition~3.2),
\begin{equation}\label{eqn2}
{\frac {HH^k_n(A/I)} {HH^k_n(A)}} \cong 
{\oplus}_{1 \le 2j \le n}H^{n-2j} \left({\frac {F^{n-j}_I({\Omega}^*_{A/k})}
{F^{n-j+1}_I({\Omega}^*_{A/k})}}\right),
\end{equation}
where $F_I{{\Omega}^*_{A/k}}$ is a filtration for the DeRham complex
${\Omega}^*_{A/k}$ whose successive quotients are given by
$${\left({\frac {F^{n-j}_I({\Omega}^*_{A/k})}
{F^{n-j+1}_I({\Omega}^*_{A/k})}} \right)}_{n-2j} = {\frac {I^j} {I^{j+1}}}
{{\otimes}_A} {\Omega}^{n-2j}_{A/k} \ {\rm and}$$   
$${\left({\frac {F^{n-j}_I({\Omega}^*_{A/k})}
{F^{n-j+1}_I({\Omega}^*_{A/k})}} \right)}_{n-2j+1} = {\frac {I^{j-1}} {I^j}}
{{\otimes}_A} {\Omega}^{n-2j+1}_{A/k}$$   
Note that since $I$ is an invertible ideal, all its powers are also
invertible, and hence ~\ref{eqn2} holds for all powers of $I$. Since the
lemma is trivial for $n = 0$, we can assume that $n$ is positive, and
so is $j$. In this case, we see that that the natural map
${\frac {{(I^2)}^j} {{(I^2)}^{j+1}}} {{\otimes}_A} {\Omega}^{n-2j}_{A/k} 
\lra {\frac {I^j} {I^{j+1}}} {{\otimes}_A} {\Omega}^{n-2j}_{A/k}$ is zero,
and hence by
comparing ~\ref{eqn2} for $I$ and $I^2$, we see that for $1 \le 2j \le n$,
the map  
$$H^{n-2j} \left({\frac {F^{n-j}_{I^2}({\Omega}^*_{A/k})}
{F^{n-j+1}_{I^2}({\Omega}^*_{A/k})}}\right) \lra 
H^{n-2j} \left({\frac {F^{n-j}_I({\Omega}^*_{A/k})}
{F^{n-j+1}_I({\Omega}^*_{A/k})}}\right)$$
is zero. Now, we use ~\ref{eqn2} again to finish the proof.
$\hfill \square$
\\
\\
Let $k$ be a field of char. 0. For any ideal $I$ of a $k$-algebra $A$,
let ${\Omega}^i_{(A, I^n)/k}$ denote the kernel of the natural surjection
${\Omega}^i_{A/k} \surj {\Omega}^i_{(A/I)/k}$.
\begin{lemma}\label{HH0}
Let $A$ be a reduced $k$-algebra, and let $B$ be the
normalisation of $A$. Let $I$ be a conducting ideal for the normalisation.
For any $i \ge 1$, the map 
$${\frac {{\Omega}^i_{(B, I^{i+1})/k}} 
{{\Omega}^i_{(A, I^{i+1})/k}}} \lra {\frac {{\Omega}^i_{(B, I)/k}} 
{{\Omega}^i_{(A, I)/k}}}$$ is zero.
\end{lemma}
{\bf Proof.}
We first observe from the universal property of the module of Kahler
differentials that ${\Omega}^i_{(A, I)/k}$ is the submodule of
${\Omega}^i_{A/k}$, generated by the exterior forms of the type
${a_0}{da_1} \wedge \cdots \wedge {da_i}$, where $a_p \in A$ for all
$p$ and $a_p \in I$ for some $p$. Let $F{{\Omega}^i_{(A, I)/k}}$ denote
the submodule of ${\Omega}^i_{(A, I)/k}$ generated by the exterior forms
of the type ${a_0}{da_1} \wedge \cdots \wedge {da_i}$, with $a_p \in I$
for all $p$. Then, it is enough to show that 
\begin{equation}\label{eqn10}
{{\rm Image}({\Omega}^i_{(B, I^{i+1})/k} \rightarrow 
{\Omega}^i_{(B, I)/k})} \subset 
{{\rm Image}(F{{\Omega}^i_{(A, I)/k}} \rightarrow {\Omega}^i_{(B, I)/k})}.   
\end{equation}
We prove this by induction on $i$. \\
For $i = 1$, let $w = adb$ with $a$ or $b$ in $I^2$. If $a \in I^2$, then
can assume $a = {a_1}{a_2}$ with $a_p \in I$. In that case, one gets
${a_1}{a_2}db = a_1(d(a_2b) - bd(a_2))$, which is clearly in
$F{{\Omega}^1_{(A, I)/k}}$. If $b \in I^2$, one proceeds
similarly. This proves $i =1$ case. Suppose now that ~\ref{eqn10} holds
for all $j \le {i-1}$ with $i > 1$. Put 
$w = {a_0}{da_1} \wedge \cdots \wedge {da_i}$ with some $a_p$ in $I^{i+1}$.
\\
{\bf Case 1.} $p = 0$ \\
Can assume $a_0 = {a^1_0} \cdots {a^{i+1}_0}$. Then
\vspace*{3mm}
$$
\begin{array}{lll}
{{a^1_0} \cdots {a^{i+1}_0}}{{da_1} \wedge \cdots \wedge {da_i}} & = &
({{a^1_0} \cdots {a^i_0}}{{da_1} \wedge \cdots \wedge {da_{i-1}}})
({a^{i+1}_0}{da_i}) \\
& = &  ({{a^1_0} \cdots {a^i_0}}{{da_1} \wedge \cdots \wedge {da_{i-1}}})
(d({a^{i+1}_0}{a_i}) - {a_i}{d{a^{i+1}_0}}) \\
& = & ({{a^1_0} \cdots {a^i_0}}{{da_1} \wedge \cdots \wedge {da_{i-1}}
\wedge d({a^{i+1}_0}{a_i})}) - \\
& & ({{a^1_0} \cdots {a^i_0}{a_i}}{{da_1} \wedge \cdots \wedge {da_{i-1}}
\wedge {d{a^{i+1}_0}}}).
\end{array}
$$
The induction hypothesis now applies. \\
{\bf Case 1.} $p > 0$. \\
The proof is exactly along the lines of case 1.
$\hfill \square$
\\
\begin{lemma}\label{HH3}
Let $A$ be a reduced ring which is essentially of finite type over a field
$K$ of char. 0, and let $B$ be the smooth normalisation of $A$. Let $I$ be
a conducting ideal for the normalisation which is invertible in $B$.
Let $k \subset K$ be a subfield. Then, for any $i \ge 1$, the natural map
$${\frac {HH^k_i(B, I^n)} {HH^k_i(A, I^n)}} \lra 
{\frac {HH^k_i(B, I)} {HH^k_i(A, I)}}$$ is zero for all sufficiently large
$n$.
\end{lemma}
{\bf Proof.}
We shall in fact show that this holds for all $n \ge {(i+1)}^2$. 
Consider the exact sequence for relative Hochschild homology
$$0 \rightarrow {\frac {HH^k_{i+1}(B/ {I^n})} {HH^k_{i+1}(B)}} 
\rightarrow HH^k_i(B, I^n) \rightarrow
{\rm Ker}(HH^k_i(B) \rightarrow HH^k_i(B/{I^n})) \rightarrow 0.$$
Since $S$ is smooth, we have seen that $HH^k_i(B) = HH^{k, (i)}_i(B)
= {\Omega}^i_{B/k}$, and hence from the naturality of Hodge decomposition
on Hochschild homology, we have 
${\rm Ker}(HH^k_i(B) \rightarrow HH^k_i(B/{I^n})) \cong
{\rm Ker}({\Omega}^i_{B/k} \rightarrow {\Omega}^i_{(B/{I^n})/k}) =
{\Omega}^i_{(B, I^n)/k}$. Thus, we get a diagram of exact sequences
$$
\xymatrix@C.6pc{
0 \ar[r] &  {\frac {HH^k_{i+1}(A/{I^n})} {HH^k_{i+1}(A)}} \ar[r] \ar[d] &
HH^k_i(A, I^n) \ar[d] \ar[r] &   
{{\rm Ker}(HH^k_i(A) \rightarrow HH^k_i(A/{I^n}))} \ar[d] \ar[r] & 0 \\
0 \ar[r] &  {\frac {HH^k_{i+1}(B/{I^n})} {HH^k_{i+1}(B)}} \ar[r] &
HH^k_i(B, I^n) \ar[r] &   
{\Omega}^i_{(B, I^n)/k} \ar[r] & 0. 
}
$$
Taking quotients, we get exact sequence
$${\frac {HH^k_{i+1}(B/{I^n})} {{HH^k_{i+1}(B)} + {HH^k_{i+1}(A/{I^n})}}}
\rightarrow {\frac {HH^k_i(B, I^n)} {HH^k_i(A, I^n)}} \rightarrow
{\frac {{\Omega}^i_{(B, I^n)/k}} 
{{\rm Ker}(HH^k_i(A) \rightarrow HH^k_i(A/{I^n}))}} \rightarrow 0.$$
Furthermore, since ${\Omega}^i_{A/k} = HH^{k, (i)}_i(A)$ ([LO]), and 
similarly for
other rings, we see from the naturality of Hodge decomposition that
${\frac {{\Omega}^i_{(B, I^n)/k}} 
{{\rm Ker}(HH^k_i(A) \rightarrow HH^k_i(A/{I^n}))}} \cong
{\frac {{\Omega}^i_{(B, I^n)/k}} 
{{\Omega}^i_{(A, I^n)/k}}}$. 
Comparing above exact sequence for $n = 1, n= {(i+1)}$ and $n = {(i+1)}^2$, 
we get a diagram
$$
\xymatrix@C.6pc{
{\frac {HH^k_{i+1}(B/{I^{{(i+1)}^2}})} {{HH^k_{i+1}(B)} + 
{HH^k_{i+1}(A/{I^{{(i+1)}^2}})}}}
\ar[r] \ar[d] & 
{\frac {HH^k_i(B, I^{{(i+1)}^2})} 
{HH^k_i(A, I^{{(i+1)}^2})}} \ar[r] \ar[d] &
{\frac {{\Omega}^i_{(B, I^{{(i+1)}^2})/k}}
{{\Omega}^i_{(A, I^{{(i+1)}^2})/k}}} \ar[r] 
\ar[d] & 0 \\
{\frac {HH^k_{i+1}(B/{I^{i+1}})} {{HH^k_{i+1}(B)} + 
{HH^k_{i+1}(A/{I^{i+1}})}}}
\ar[r] \ar[d] & 
{\frac {HH^k_i(B, I^{i+1})} 
{HH^k_i(A, I^{i+1})}} \ar[r] \ar[d] &
{\frac {{\Omega}^i_{(B, I^{i+1})/k}} 
{{\Omega}^i_{(A, I^{i+1})/k}}} \ar[r] \ar[d] & 0 \\
{\frac {HH^k_{i+1}(B/I)} {{HH^k_{i+1}(B)} + 
{HH^k_{i+1}(A/I)}}}
\ar[r] & {\frac {HH^k_i(B, I)} 
{HH^k_i(A, I)}} \ar[r] &
{\frac {{\Omega}^i_{(B, I)/k}} 
{{\Omega}^i_{(A, I)/k}}}  \ar[r] & 0 
}
$$
The two vertical maps on the left are zero by lemma~\ref{HH2}, and the two
vertical maps on the right are zero by lemma~\ref{HH0}. A diagram chase
shows that the composite map in the middle is zero.
$\hfill \square$
\\
\begin{lemma}\label{HH4}
Under the conditions of lemma~\ref{HH3}, the map
$$HH^k_1(A, I^n) \lra HH^k_1(B, I^n)$$ is injective for all sufficiently
large $n$.
\end{lemma}
{\bf Proof.}
For any subring $R \inj A$, let $D^R_*(A, I)$ be the relative Andre-Quillen
homology defined as the homology groups of the complex
${\rm Ker}({\mathbb L}_{A/R} \surj {\mathbb L}_{(A/I)/R})$. These groups
fit into long exact sequence of relative Andre-Quillen homology.
As in [LO], there are natural maps $D^R_i(A, I) \lra HH^{R, (1)}_{i+1}
(A, I)$. Comparing thse groups using long exact sequences of Andre-Quillen
homology and Hochschild homology, and using the isomorphism
$D^R_i(A) \cong HH^{R, (1)}_{i+1}(A)$, one gets isomorphism
$D^R_i(A, I) \cong HH^{R, (1)}_{i+1}(A, I)$ for all $i \ge 0$.
Thus, we need to show that the natural map $D^k_0(A, I^n) \lra
D^k_0(B, I^n)$ is injective for all large $n$. Using the base change
long exact sequence of Andre-Quillen homology ([LO]), one gets exact
sequence
$${\Omega}^1_{K/k} \otimes {I^n} \lra D^k_0(A, I^n) \lra D^K_0(A, I^n)
\lra 0$$  
Comparing this exact sequence for $A$ and $B$, we have a commutative
diagram
$$
\xymatrix@C.6pc{
& {\Omega}^1_{K/k} \otimes {I^n} \ar[r] \ar@{=}[d] & 
D^k_0(A, I^n) \ar[d] \ar[r] & D^K_0(A, I^n) \ar[r] \ar[d] & 0 \\
& {\Omega}^1_{K/k} \otimes {I^n} \ar[r] \ar@{^{(}->}[d] & 
D^k_0(B, I^n) \ar[d] \ar[r] & D^K_0(B, I^n) \ar[r] \ar[d] & 0 \\
0 \ar[r] & {\Omega}^1_{K/k} \otimes B \ar[r] &
{\Omega}^1_{B/k} \ar[r] & {\Omega}^1_{B/K} \ar[r] & 0,
}
$$
in which all the rows are exact, and the second diagram is a part of
long exact sequence of relative Andre-Quillen homology of $B$ and $I^n$,
and using the isomorphism $D^k_0(B) \cong {\Omega}^1_{B/k}$. The bottom
sequence is exact on the left since $B$ is smooth over $K$. 
A diagram chase now shows that all the rows are exact on the left. Now,
lemma~\ref{HH1} and Snake lemma complete the proof.
$\hfill \square$
\\
\begin{lemma}\label{HH5}
Under the conditions of lemma~\ref{HH3}, the natural map
$$HH^k_1(A, B, I^n) \lra HH^k_1(A, B, I)$$ of double relative Hochschild
homology groups is zero for all sufficiently large $n$.
\end{lemma}
{\bf Proof.}
The long exact sequence of relative Hochschild homology gives exact 
sequence
$$0 \rightarrow {\frac {HH^k_2(B, I^n)} {HH^k_2(A, I^n)}} 
\rightarrow HH^k_1(A, B, I^n)
\rightarrow {\rm Ker}(HH^k_1(A, I^n) \rightarrow HH^k_1(B, I^n)) 
\rightarrow 0.$$
But the last group vanishes by lemma~\ref{HH4}. Furthermore, the map
${\frac {HH^k_2(B, I^n)} {HH^k_2(A, I^n)}} \lra 
{\frac {HH^k_2(B, I)} {HH^k_2(A, I)}}$ is zero for all large $n$ by 
lemma~\ref{HH3}. 
$\hfill \square$
\\
\begin{corollary}\label{HC1}
Under the conditions of lemma~\ref{HH3}, the natural map
$$HC^k_1(A, B, I^n) \lra HC^k_1(A, B, I)$$ of double relative cyclic homology
groups is zero for all sufficiently large $n$.
\end{corollary}
{\bf Proof.}
In view of the above lemma, it's enough to show that the natural map
$HH^k_1(A, B, I^n) \lra HC^k_1(A, B, I^n)$ is surjective for all $n$.
But, the SBI sequence ([LO]) of double relative Hochschild and cyclic
homology groups gives exact sequence
$$HH^k_1(A, B, I^n) \lra HC^k_1(A, B, I^n) \stackrel {S} 
\lra HC^k_{-1}(A, B, I^n).$$ 
Another exact sequence of relative cyclic homology gives
exact sequence
$$HC^k_0(A, I^n) \lra HC^k_0(B, I^n) \lra  
HC^k_{-1}(A, B, I^n) \lra HC^k_{-1}(A, I^n).$$
However, $HC^k_0(A, I^n) \cong I^n \cong HC^k_0(B, I^n)$, and
$HC^k_{-1}(A, I^n) = 0$.
This finishes the proof.
$\hfill \square$
\\
\\
{\bf Remark.} We point out here that it is already known that the map
$HC^k_0(A, B, I^2) \rightarrow HC^k_0(A, B, I)$ is zero ([CGW]).
\\

We conclude this section with the following two lemmas.
\begin{lemma}\label{TH}
Let $A$ be an one-dimensional local domain, essentially of finite type
over an algebraically closed field of char. 0, let $B$ be the normalisation
of $A$. Let $k \subset K$ be any subfield. Then the natural map
$${\rm Ker}(HH^k_1(A) \rightarrow HH^k_1(B)) \lra 
{\rm Ker}(HH^K_1(A) \rightarrow HH^K_1(B))$$
is an isomorphism.
\end{lemma}
{\bf Proof.}
We first observe that $HH^k_1(A) = {\Omega}^1_{A/k}$ and similarly for $B$.
Thus, we can work with Kahler differentials. The above map is already
injective, so we need to show only surjectivity.
Note that since $B$ is regular, ${\Omega}^1_{B/k}$ is a free $B$-module
(not necessarily finitely generated). Thus, the map
${\Omega}^1_{K/k} {{\otimes}_K} A \inj  
{\Omega}^1_{K/k} {{\otimes}_K} B \rightarrow {\Omega}^1_{B/k}$ is injective.
Furthermore, since $K$ is algebraically closed, the map 
${\Omega}^1_{K/k} \rightarrow {\Omega}^1_{A/k} \rightarrow
{\Omega}^1_{K/k}$ is naturally split. In particular, 
${\Omega}^1_{K/k}$ is naturally a direct summand of ${\Omega}^1_{A/k}$.
\\
Put ${\ov {{\Omega}^1_{A/k}}} = {\rm Ker}({\Omega}^1_{A/k} \rightarrow
{\Omega}^1_{B/k})$. We define $\ov {{\Omega}^1_{A/K}}$ similarly.
Let $F$ denote the field of fractionas of $A$. Then it is easy to see
that ${\ov {{\Omega}^1_{A/k}}} = {\rm Tor}^1_A(F/A, {\Omega}^1_{A/k})$,
and one has similar interpretation for $\ov {{\Omega}^1_{A/K}}$.
This follows because the map ${\Omega}^1_{B/k} \lra {\Omega}^1_{F/k}$
is injective. Thus, we need to show that the map
${\rm Tor}^1_A(F/A, {\Omega}^1_{A/k}) \lra 
{\rm Tor}^1_A(F/A, {\Omega}^1_{A/K})$ is surjective. But, using the 
exact sequence
$$0 \lra {\Omega}^1_{K/k} {{\otimes}_K} A \lra {\Omega}^1_{A/k} \lra
{\Omega}^1_{A/K} \lra 0,$$
one gets a long exact sequence
$${\rm Tor}^1_A(F/A, {\Omega}^1_{A/k}) \rightarrow 
{\rm Tor}^1_A(F/A, {\Omega}^1_{A/K}) \rightarrow
{\Omega}^1_{K/k} {{\otimes}_K} F/A \stackrel {\phi} \rightarrow 
{\Omega}^1_{A/k} {{\otimes}_A} F/A.$$
Hence, it is enough to show that $\phi$ is injective.
However, one has a factorisation ${\Omega}^1_{K/k} {{\otimes}_K} F/A 
\rightarrow {\Omega}^1_{A/k} {{\otimes}_A} F/A \rightarrow
A {{\otimes}_K} ({\Omega}^1_{A/k} {{\otimes}_A} F/A)$. Thus, it is enough
to show that the composite map is injective.
However,  
$${\Omega}^1_{K/k} {{\otimes}_K} F/A =
({\Omega}^1_{K/k} {{\otimes}_K} A) {{\otimes}_A} F/A 
\inj ({\Omega}^1_{A/k} {{\otimes}_K} A) {{\otimes}_A} F/A$$
$$\ \ \ \ \ \ \ \ \ \ \ 
\cong A {{\otimes}_K} ({\Omega}^1_{A/k} {{\otimes}_A} F/A).$$
Here, the injective arrow follows because ${\Omega}^1_{K/k}$
is naturally a direct summand of ${\Omega}^1_{A/k}$ as observed
before. This proves the desired injectivity.
$\hfill \square$
\\
\begin{lemma}\label{TH1}
Let $A$ and $B$ be as in lemma~\ref{TH}. Let $\mathfrak m$ and
$\mathfrak M$ denote the Jacobson radicals of $A$ and $B$ respectively.
Then the natural map
$${\rm Ker}(HH^k_1(A, {\mathfrak m}) \rightarrow HH^k_1(B, {\mathfrak M}))
\lra {\rm Ker}(HC^k_1(A, {\mathfrak m}) \rightarrow 
HC^k_1(B, {\mathfrak M})).$$
is injective.
\end{lemma}
{\bf Proof.}
We consider the following diagram of exact sequences coming from 
the SBI-sequence
$$
\xymatrix@C.6pc{
0 \ar[r] & HC^k_0(A, {\mathfrak m}) \ar[d] \ar[r] &
HH^k_1(A, {\mathfrak m}) \ar[d] \ar[r] &
HC^k_1(A, {\mathfrak m}) \ar[d] \ar[r] & 0 \\
0 \ar[r] & HC^k_0(B, {\mathfrak M}) \ar[r] &
HH^k_1(B, {\mathfrak M}) \ar[r] &
HC^k_1(B, {\mathfrak M}) \ar[r] & 0,
}
$$
where the first arrow from left in the bottom sequence is injective because
the composite map ${\mathfrak M} = HC^k_0(B, {\mathfrak M}) \lra
HH^k_1(B, {\mathfrak M}) \lra {\Omega}^1_{B/k}$ is injective as $B$ is
regular. Also, $HC^k_0(A, {\mathfrak m}) \cong {\mathfrak m}$, again
using the long exact sequence for relative cyclic homology. Thus the 
left-most vertical map is just the inclusion ${\mathfrak m} 
\inj {\mathfrak M}$. A diagram chase now proves the lemma.
$\hfill \square$
\\    
\section{Mayer-Vietoris sequences in $K$-theory and Cyclic homology}
Our goal in this section is to establish some Mayer-Vietoris type exact 
sequences in $K$-theory and cyclic homology. These sequences will be one of
our main tools to prove main theorems. Let $A$ be an 
one-dimensional reduced local ring, which is essentially of finite type
over an algebraically closed field $K$ of characteristic zero. Let 
$\mathfrak m$ denote the Jacobson radical of $A$. Let $B$ be the normalisation
of $A$ with the Jacobson radical $\mathfrak M$. Note that $B$ is a direct
product of regular semi-local domains. Then, for any radical ideal $I$ of $B$,
one has a fibration of $K$-theory spectra
$$K(B, I) \lra K(B, {\mathfrak M}) \lra K(B/I, {\mathfrak M}/I).$$
\begin{lemma}\label{KMV}
Let $A$ be a reduced local ring as above with the maximal ideal 
${\mathfrak m}$, and let $B$ be the normalisation of $A$ with Jacobson 
radical $\mathfrak M$.
Then for any conducting ideal $I$, one has `Mayer-Vietoris' exact
sequences
$$K_2(A, {\mathfrak m}) \lra K_2(B, {\mathfrak M}) \oplus 
K_2(A/I, {\mathfrak m}/I)  \lra K_2(B/I, {\mathfrak M}/I) \lra 0.$$
$$HC_1(A, {\mathfrak m}) \lra HC_1(B, {\mathfrak M}) \oplus 
HC_1(A/I, {\mathfrak m}/I)  \lra HC_1(B/I, {\mathfrak M}/I) \lra 0.$$
\end{lemma}
{\bf Proof.}
From the above fibration of $K$-theory spectra, one has diagrams of exact
sequences
$$
\xymatrix@C.4pc{
& K_2(A, I) \ar[d] \ar[r] & K_2(A, {\mathfrak m}) \ar[d] \ar[r] &
K_2(A/I, {\mathfrak m}/I) \ar[d] \ar[r] & 0 \\
K_3(B/I, {\mathfrak M}/I) \ar[r] \ar[dr] & K_2(B, I) \ar[r] \ar[d] & 
K_2(B, {\mathfrak M}) \ar[r]  & K_2(B/I, {\mathfrak M}/I) \ar[r] & 
0 \\
& K_1(A, B, I) & & & 
}
$$
$$
\xymatrix@C.4pc{
& HC_1(A, I) \ar[d] \ar[r] & HC_1(A, {\mathfrak m}) \ar[d] \ar[r] &
HC_1(A/I, {\mathfrak m}/I) \ar[d] \ar[r] & 0 \\
HC_2(B/I, {\mathfrak M}/I) \ar[r] \ar[dr] & HC_1(B, I) \ar[r] \ar[d] & 
HC_1(B, {\mathfrak M}) \ar[r]  & HC_1(B/I, {\mathfrak M}/I) \ar[r] & 
0 \\
& HC_0(A, B, I) & & & 
}
$$
$${\rm Diagram} \ 2$$
The surjectivity of the last horizontal map in the first diagram folllows 
since the map $B/I \lra B/{{\mathfrak M}}$ is split surjective, and then 
compare long exact $K$-theory sequence for pairs $(B, {\mathfrak M})$ and
$(B/I, {\mathfrak M}/I)$. Similar argument holds in the second diagram.
This also proves the surjectivity of the last maps
in the lemma. Now a diagram chase shows that it is enough to prove that
the slanted arrows in both diagrams are surjective. However,
we know that by [GO] and [CO], there are isomorphisms
$$K_3(B/I, {\mathfrak M}/I) \cong HC_2(B/I, {\mathfrak M}/I),
\ {\rm and} \ K_1(A, B, I) \cong HC_0(A, B, I).$$
Here, all Hochschild and Cyclic homologies are taken with respect to the
field of rational numbers $\Q$. Furthermore, one has a commutative diagram
$$
\xymatrix{
HH_2(B/I, {\mathfrak M}/I) \ar[r] \ar[d] & HH_0(A, B, I) \ar@{=}[d] \\
HC_2(B/I, {\mathfrak M}/I) \ar[r] & HC_0(A, B, I),
}
$$
where the right eqaulity follows from the SBI-sequence of double relative
Hochschild and Cyclic homology ([LO]) or by direct computation.
But, from the proof of Theorem~1.2 of [CGW], the map $HH_2(B/I, 
{\mathfrak M}/I) \lra HH_0(A, B, I)$ is surjective. This proves desired
surjectivity of both the slanted arrows.
$\hfill \square$
\\
\begin{corollary}
With notations as in the above lemma, the maps
$${\rm Ker}(K_2(A, {\mathfrak m}) \rightarrow K_2(B, {\mathfrak M})) 
\lra {\rm Ker}(K_2(A/I, {{\mathfrak m}/I}) \rightarrow K_2(B/I, 
{{\mathfrak M}/I})),$$
$${\rm Ker}(HC_1(A, {\mathfrak m}) \rightarrow HC_1(B, {\mathfrak M})) 
\lra {\rm Ker}(HC_1(A/I, {{\mathfrak m}/I}) \rightarrow HC_1(B/I, 
{{\mathfrak M}/I}))$$
are surjective.
\end{corollary}
{\bf Proof.}
Follows directly from the `Mayer-Vietoris' sequence of the lemma.
$\hfill \square$
\\
\\
{\bf Remark.} The second part of the corollary was also established in [CGW].
\\
\begin{lemma}\label{KMV1}
Under the hypothesis of lemma~\ref{KMV}, the map
$${\frac{K_3(B/{I^2}, {{\mathfrak M}/{I^2}})} {K_3(B, {\mathfrak M})}}
\lra {\frac{K_3(B/I, {{\mathfrak M}/I})} {K_3(B, {\mathfrak M})}}$$
is zero.
\end{lemma}
{\bf Proof.}
Note that $K_3(B/I, {{\mathfrak M}/I}) = K^{{\rm nil}}_3(B/I, 
{{\mathfrak M}/I})$, and the latter is a $\Q$-vector space. Hence, both
groups above remain unchanged even after we mod out torsion part
of $K_3(B, {\mathfrak M})$. Using
Adam's operations on rational relative $K$-theory as in [L] (see also
[C]), one has eigenspace decomposition
$$K_3(B/{I^n}, {{\mathfrak M}/{I^n}}) = K^{(2)}_3(B/{I^n}, 
{{\mathfrak M}/{I^n}}) \oplus K^{(3)}_3(B/{I^n}, {{\mathfrak M}/{I^n}}).$$
Further, $K^{(3)}_3(B/{I^n}, {{\mathfrak M}/{I^n}}) =
K^{M}_3(B/{I^n}, {{\mathfrak M}/{I^n}})$, where the latter is the relative
Milnor $K$-group as defined by Levine ([L]). By naturality of eigenspace
decomposition, one gets
$$
{\frac{K_3(B/{I^n}, {{\mathfrak M}/{I^n}})} {K_3(B, {\mathfrak M})}}
\cong 
{\frac{K^{(2)}_3(B/{I^n}, {{\mathfrak M}/{I^n}})} 
{K^{(2)}_3(B, {\mathfrak M})}} \oplus 
{\frac{K^{M}_3(B/{I^n}, {{\mathfrak M}/{I^n}})} 
{K^{(3)}_3(B, {\mathfrak M})}}.
$$
Now, $B$ is a direct product of regular semi-local domains in which 
all height 1 prime ideals are
pricipal, and since $\mathfrak M$ is the product of all maximal ideals,
we see the pair
$(B, {\mathfrak M})$ satisfies the MV-Property of Levine.
We 
\\
{\bf Claim.} There is a surjection
$$K^{M}_3(B, {\mathfrak M}) \surj {\rm Ker}(K^M_3(B) \surj 
K^M_3(B/{\mathfrak M})).$$
For this, we use the eigen pieces of the long exact rel. $K$-theory sequence,
to get exact sequence 
$$K^{(3)}_3(B, {\mathfrak M}) \lra {\rm Ker}(K^{(3)}_3(B) \lra 
K^{(3)}_3(B/{\mathfrak M})).$$ But all these are corresponding Milnor
$K$-groups by [L]. This proves the claim. We point out here that
the isomorphism $K^{M}_3(B, {\mathfrak M}) \cong F^3K_3(B, {\mathfrak M})
\cong K^{(3)}_3(B, {\mathfrak M})$ is known only after we mod out
two torsion elements. But as remarked in the beginning of the proof
of the lemma, this does not affect the statement of the lemma.

Now, using this claim and the the fact that the surjection $B/{I^n} \surj
B/{\mathfrak M}$ splits, one has a diagram of exact sequences
$$
\xymatrix{
& K^{M}_3(B, {\mathfrak M}) \ar[d] \ar[r]  & K^M_3(B) \ar[d] \ar[r] &
K^M_3(B/{\mathfrak M}) \ar[d] \ar[r] & 0 \\
0 \ar[r] & K^M_3(B/{I^n}, {{\mathfrak M}/{I^n}}) \ar[r] &
K^M_3(B/{I^n}) \ar[r] & K^M_3(B/{\mathfrak M}) \ar[r] & 0,
}
$$
which in turn gives a surjection
$$ K^{M}_3(B, {\mathfrak M}) \surj  K^M_3(B/{I^n}, {{\mathfrak M}/{I^n}}).$$
Applying this in the eigenspace decomposition above, we obtain
\vspace*{3mm}
$$
\begin{array}{lll}
{\frac{K_3(B/{I^n}, {{\mathfrak M}/{I^n}})} {K_3(B, {\mathfrak M})}} & =
& {\frac{K^{(2)}_3(B/{I^n}, {{\mathfrak M}/{I^n}})} 
{K^{(2)}_3(B, {\mathfrak M})}} \\
& = & {\frac{HC^{(1)}_2(B/{I^n}, {{\mathfrak M}/{I^n}})} 
{K^{(2)}_3(B, {\mathfrak M})}} \\
& = & {\frac{{\prod}_i {HC^{(1)}_2({B_{{\mathfrak M}_i}}/{I^n}, 
{{{\mathfrak M}_i}/{I^n}})}} {K^{(2)}_3(B, {\mathfrak M})}} \\
& = &  {\frac{{\prod}_i{D^{\Q}_1({\Q}[t^i]/{t^{r_i}_i}, (t_i)) \otimes
HH^{\Q}_0(k)}} {K^{(2)}_3(B, {\mathfrak M})}},
\end{array}
$$
where the sum is taken over maximal ideals of $B$.
Here, $D^{\Q}(A)$ denotes the Andre-Quillen homology of a $\Q$-algebra $A$
([LO])
and the last equality follows from the computation of the
Cyclic homology of truncated polynomial algebras in ([LO], sec. 4.6).  
Now, the proof of the lemma follows from the following
\begin{sublemma}\label{Andre}
Let $B_r$ denote the truncated polynomial ring ${{\Q}[t]}/(t^r)$. Then, the
map $D^{\Q}_1(B_{2r}) \lra D^{\Q}_1(B_r)$ is zero.
\end{sublemma}
{\bf Proof.}
Since we are dealing with rational coefficients, we shall
ignore the index $\Q$ in this proof. Note that $D_1(B_r) = 
{\frac{D_1(B_r)} {D_1({\Q}[t])}}$ (${\Q}[t]$ is smooth over $\Q$),
which in turn is same as 
${\frac{HH^{(1)}_2(B_r)} {HH^{(1)}_2({\Q}[t])}}$ by [LO]. But the map
$${\frac{HH_2(B_{2r})} {HH_2({\Q}[t])}} \lra
{\frac{HH_2(B_r)} {HH_2({\Q}[t])}}$$ is zero by lemma~\ref{HH2}.
$\hfill \square$
\begin{corollary}\label{CKMV1}
Let $F_2(B, I)$ denote Ker$(K_2(B, I) \lra K_2(B, {\mathfrak M}))$. Then 
the map $F_2(B, I^2) \lra F_2(B, I)$ is zero.
\end{corollary}
{\bf Proof.}
This follows once we observe that 
$$F_2(B, I) = {\frac{K_3(B/I, {{\mathfrak M}/I})} {K_3(B, {\mathfrak M})}},$$
using the $K$-theory long exact sequence for the map of pairs
$(B, I) \lra (B, {\mathfrak M})$, and then using above lemma.
$\hfill \square$
\\

The following is our main result of this section, which is a stronger 
version of lemma~\ref{KMV}. 
\begin{theorem}\label{KMV2}
Consider the hypothesis of lemma~\ref{KMV}. Then, there exists a conducting
ideal $I$ such that one has `Mayer-Vietoris' exact sequences
\begin{equation}
0 \rightarrow  K_2(A, {\mathfrak m}) \rightarrow  K_2(B, {\mathfrak M}) 
\oplus K_2(A/I, {\mathfrak m}/I)  \rightarrow  K_2(B/I, {\mathfrak M}/I) 
\rightarrow 0.
\end{equation}
\begin{equation}
0 \rightarrow HC_1(A, {\mathfrak m}) \rightarrow HC_1(B, {\mathfrak M}) 
\oplus HC_1(A/I, {\mathfrak m}/I)  \rightarrow HC_1(B/I, {\mathfrak M}/I) 
\rightarrow 0.
\end{equation}
\end{theorem}
{\bf Proof.}
In view of lemma~\ref{KMV}, we only need to prove injectivity of first
maps in both sequences for some conducting ideal $I$. We shall in fact
show that given a conducting ideal $I$, this holds for all sufficiently
large powers of $I$. We fix some notations before beginning the proof.
For any conducting ideal $I$, let 
$$F(I) : = {{\rm Ker}(K_2(A, {\mathfrak m}) \rightarrow  
K_2(A/I, {\mathfrak m}/I))} \bigcap 
{{\rm Ker}(K_2(A, {\mathfrak m}) \rightarrow  
K_2(B, {\mathfrak M}))},$$
$$F(A) : = {{\rm Ker}(K_2(A, {\mathfrak m}) \rightarrow  
K_2(B, {\mathfrak M}))}, \ {\rm and}$$
$$F(A/I) : =  {{\rm Ker}(K_2(A/I, {\mathfrak m}/I) \rightarrow  
K_2(B/I, {\mathfrak M}/I))}.$$
Then, lemma~\ref{KMV} implies that one has a short exact sequence
\begin{equation}\label{exact1}
0 \lra F(I) \lra F(A) \lra F(A/I) \lra 0.
\end{equation}
We consider the diagram of exact sequences
$$
\xymatrix{
& & K_2(A, B, I) \ar[d] & \\
0 \ar[r] & F_2(A, I) \ar[r] \ar[d] & K_2(A, I) \ar[d] \ar[r] & 
K_2(A, {\mathfrak m}) \ar[d] \\
0 \ar[r] & F_2(B, I) \ar[r] \ar[dr] & K_2(B, I) \ar[r] \ar@{->>}[d] & 
K_2(B, {\mathfrak M}) \\
& & K_1(A, B, I) & 
}
$$ 
where the groups on the left are as defined in corollary~\ref{CKMV1}. 
First we claim that the slanted arrow in this diagram is surjective.
But this follows directly once we chase diagram~2 and observe in the
proof of lemma~\ref{KMV} that the slanted arrow in that diagram is 
surjective. Thus, a diagram chase above gives exact sequence
$$K_2(A, B, I) \longrightarrow F(I) \longrightarrow 
\frac{F_2(B, I)} {F_2(A, I)} \stackrel{{\beta}_I} 
\longrightarrow  K_1(A, B, I) \rightarrow 0 $$
However, we have natural isomorphism $K_2(A, B, I) \cong HC_1(A, B, I)$
by Cortinas' theorem ([CO]). Using this in this exact sequence, and comparing
the resulting sequences for various powers of $I$, we get a diagram
$$
\xymatrix{
HC_1(A, B, I^n) \ar[d] \ar[r] & F(I^n) \ar[d] \ar[r] & 
\frac{F_2(B, I^n)} {F_2(A, I^n)} \ar[d] \\
HC_1(A, B, I) \ar[r] & F(I) \ar[r] & 
\frac{F_2(B, I)} {F_2(A, I)}
}
$$
By corollary~\ref{CKMV1}, the right vertical map is zero for $n \ge 2$, 
and the left verical map is zero for $n > > 0$ by corollary~\ref{HC1}.
Let $N$ be an integer such that both these maps are zero. Now, we repeat 
the same argument in the above diagram with $I$ replaced by $J = I^N$ 
to get a diagram
$$
\xymatrix{
HC_1(A, B, J^n) \ar[d] \ar[r] & F(J^n) \ar[d] \ar[r] & 
\frac{F_2(B, J^n)} {F_2(A, J^n)} \ar[d] \\
HC_1(A, B, J) \ar[r] \ar[d] & F(J) \ar[r] \ar[d] & 
\frac{F_2(B, J)} {F_2(A, J)} \ar[d] \\
HC_1(A, B, I) \ar[r] & F(I) \ar[r] & 
\frac{F_2(B, I)} {F_2(A, I)},
}
$$   
such that all the vertical maps on the left and right ends are zero. A diagram
chase above now shows that for all $n > > 0$, the map
$F(I^n) \lra F(I)$ is zero. Applying this in ~\ref{exact1}, we see that
$F(A) \lra F(A/I^n)$ is an isomorphism for all large powers of $I$. This
proves the exactness of first sequence. The case of cyclic homology follows
along exactly similar lines. In fact, we have reduced $K$-theory problem to 
Cyclic homology problem in the above proof.
$\hfill \square$
\\
\begin{corollary}\label{MVS}
Let $(A, {\mathfrak m})$ be a reduced one-dimensional local ring, and let
$B$ be a reduced semi-local ring containing $A$ and contained in the 
normalisation of $A$. Let $\mathfrak M$ be the Jacobson radical of $B$. Then
there is a conducting ideal $I \subset B$ such that one has Mayer-Vietoris 
exact sequences as in Theorem~\ref{KMV2}.
\end{corollary}
{\bf Proof.}
To prove the injectivity of the map $K_2(A, {\mathfrak m}) \lra  
K_2(B, {\mathfrak M}) \oplus K_2(A/I, {\mathfrak m}/I)$, observe
that we can choose a $I$ to be a conducting ideal for the normalisation
of $A$ and hence it will also be conducting ideal for $A \lra B$. This
reduces the proof to the case when $B$ is the normalisation of $A$. 
To prove the exactness of the sequence
$$K_2(A, {\mathfrak m}) \lra  
K_2(B, {\mathfrak M}) \oplus K_2(A/I, {\mathfrak m}/I) \lra 
K_2(B/I, {\mathfrak M}/I) \lra 0,$$
we use the exact sequence (which always holds) from the Diagram~3
$$K_1(A, B, I) \lra {\frac{K_2(B, {\mathfrak M})} {K_2(A, {\mathfrak m})}}
\lra {\frac{K_2(B/I, {\mathfrak M}/I)} {K_2(A/I, {\mathfrak m}/I)}}
\lra 0$$
and use the Cortinas' isomorphism $K_1(A, B, I) \cong HC_0(A, B, I) 
\cong {I/{I^2}} \otimes {\Omega}^1_{B/A}$, which holds even if $B$ is not
normal. Now, we copmare this exact sequence for $I$ and $I^2$ and
argue as before to finish the proof. The case of Cyclic homology is along
the similar lines.
$\hfill \square$
\\
\section{Proofs of main theorems}
{\bf {\large {Proof of Theorem~\ref{thm1}.}}}
Let $A$ be an one-dimensional local domain, essentially of finite type over
an algebraically closed field $K$ of characteristic 0. Let $\mathfrak m$
denote the Jacobson radical of $A$. Let $B$ be the normalisation of $A$
with Jacobson radical $\mathfrak M$. Let $F$ be the field of fractions of 
$A$. Assume that the `Artinian Geller Conjecture' holds, and $A$ is singular. 
Since $B$ is a regular semi-local domain, the map $K_2(B) \lra K_2(F)$ is
injective by Quillen's proof of Gersten conjecture. Hence it suffices to 
prove that the map $K_2(A) \lra K_2(B)$ is not injective. Consider the 
commutative diagram with exact rows
$$
\xymatrix@C.6pc{
0 \ar[r] & K_2(A, {\mathfrak m}) \ar[r] \ar[d] & K_2(A) \ar[r] \ar[d] &
K_2(A/{\mathfrak m}) \ar[r] \ar[d] & 0 \\
& K_2(B, {\mathfrak M}) \ar[r] & K_2(A) \ar[r] &
K_2(B/{\mathfrak M}) \ar[r] & 0. 
}
$$
$${\rm Diagram} \ 3$$
Here, the map $K_2(A, {\mathfrak m}) \lra K_2(A)$ is injective since
$K \lra A \lra A/{\mathfrak m} = K$ is split, as $K$ is algebraically closed.
By the same reason, the right-most vertical map is injective, since
$B/{\mathfrak M}$ is some copies of $K$. A diagram chase shows that it is
enough to show that the left-most vertical map is not injective.
We choose a conducting ideal $I$ for the normalisation such that 
$I \subset {\mathfrak m}^2$ and moreover, one has the 
Mayer-Vietoris exact sequences as in Theorem~\ref{KMV2}. Thus, we have
\begin{equation}\label{KH1}
{\rm Ker}(K_2(A, {\mathfrak m}) \rightarrow K_2(B, {\mathfrak M})) 
\cong {\rm Ker}(K_2(A/I, {{\mathfrak m}/I}) \rightarrow K_2(B/I, 
{{\mathfrak M}/I})).
\end{equation}
Using the same diagram as above with $A$ (resp. $B$) replaced with $A/I$
(resp. $B/I$), and obsering that $K_3(B/I) \surj K_3(B/{\mathfrak M})$, 
we see that
\begin{equation}\label{KH2}
{\rm Ker}(K_2(A/I) \rightarrow K_2(B/I)) = 
{\rm Ker}(K_2(A/I, {\mathfrak m}/I) 
\rightarrow K_2(B/I, {\mathfrak M}/I)).
\end{equation} 
Now, if $A$ is singular, then $\mathfrak m$ is not a principal ideal, 
and since $I \subset {\mathfrak m}^2$, Nakayama's Lemma implies that
${\mathfrak m}/I$ is also not a principal ideal in $A/I$. In particular,
$A/I$ is not a principal ideal algebra though it is a subalgebra of $B/I$,
which is a principal ideal algebra. Hence by `Artinian Geller Conjecture',
${\rm Ker}(K_2(A/I) \rightarrow K_2(B/I)) \neq 0$. Now, we use ~\ref{KH2}
and then ~\ref{KH1} to finish the proof.
$\hfill \square$
\\  
\\
{\bf {\large {Proof of Theorem~\ref{thm2}.}}}
We first observe that for any subfield $k \subset K$, 
$HH^k_2(B/{\mathfrak M}) =  HH^{k, (1)}_2(B/{\mathfrak M}) \oplus
HH^{k, (2)}_2(B/{\mathfrak M})$ by the Hodge decomposition on Hochschild
homology. But $HH^{k, (1)}_2(B/{\mathfrak M}) = D^k_1(B/{\mathfrak M})
= 0$ ([C]). Also, $HH^{k, (2)}_2(B/{\mathfrak M}) =
{\Omega}^2_{({B/{\mathfrak M}})/k}$, and
${\Omega}^2_{B/k} \surj {\Omega}^2_{({B/{\mathfrak M}})/k}$. In particular,
$HH^k_2(B) \surj HH^k_2(B/{\mathfrak M})$.
Thus, using the long exact sequence for relative Hochschild homology,
we see that $${\rm Ker}({\Omega}^1_{A/k} \lra {\Omega}^1_{B/k}) =
{\rm Ker}(HH^k_1(A) \lra HH^k_1(B)) \cong $$
$$\ \ \ \ \ \ \ \ \ \ \ \ \ \ \ \ \ \ \ \ \ \ 
{\rm Ker}(HH^k_1(A, {\mathfrak m}) \lra HH^k_1(B, {\mathfrak M})).$$  
Now, suppose that Berger's conjecture holds, and $A$ is singular. Then
the map ${\Omega}^1_{A/K} \lra {\Omega}^1_{B/K}$ is not injective.
Hence, by lemma~\ref{TH}, the map $HH^{\Q}_1(A) \rightarrow
HH^{\Q}_1(B)$ is not
injective, and the above isomorphism implies that the map
$HH^{\Q}_1(A, {\mathfrak m}) \rightarrow HH^{\Q}_1(B, {\mathfrak M}))$
is not injective.
Now, we use lemma~\ref{TH1} to conclude that
\begin{equation}\label{KH3}
{\rm Ker}(HC^{\Q}_1(A, {\mathfrak m}) \lra 
HC^{\Q}_1(B, {\mathfrak M})) \neq 0.
\end{equation}

We choose a conducting ideal $I$ for the normalisation of $A$ so that we
have Mayer-Vietoris exact sequences as in Theorem~\ref{KMV2}.  
Then we get an isomorphism as in ~\ref{KH1} and also an isomorphism
\begin{equation}\label{KH3}
{\rm Ker}(HC^{\Q}_1(A, {\mathfrak m}) \rightarrow 
HC^{\Q}_1(B, {\mathfrak M})) \cong 
{\rm Ker}(HC^{\Q}_1(A/I, {{\mathfrak m}/I}) \rightarrow 
HC^{\Q}_1(B/I, {{\mathfrak M}/I})).
\end{equation}  
However, by Goodwillie's theorem ([GO]), the maps
$$K_2(A/I, {{\mathfrak m}/I}) \lra HC^{\Q}_1(A/I, {{\mathfrak m}/I}),
\ {\rm and}$$
$$K_2(B/I, {{\mathfrak M}/I})) \lra  HC^{\Q}_1(B/I, {{\mathfrak M}/I})$$
are isomorphisms.
Now, we combine ~\ref{KH3}, ~\ref{KH2} and ~\ref{KH1} to conclude that
${\rm Ker}(K_2(A, {\mathfrak m}) \rightarrow K_2(B, {\mathfrak M}))$
is not zero. But we have seen in diagram~3 that this group injects
inside ${\rm Ker}(K_2(A) \rightarrow K_2(B))$. This proves the theorem.
$\hfill \square$
\\   
\\
{\bf {\large {Proof of Corollary~\ref{cor}.}}}
The corollary follows from Theorem~\ref{thm2} since Berger's conjecture
has been verified in these cases. For example, (i) and (ii) are verified
in [CGW], (iii) in [B1], (iv) in [B2], (v) in [S], (vi) in [BA], (vii)
in [GU], and (viii) in [HW].
$\hfill \square$
\\

Department of Mathematics \\
University of California at Los Angeles \\
405 Hilgard Avenue \\
Los Angeles, CA, 90095, USA. \\
email : amalendu@math.ucla.edu \\
Tel: 1-310-825-4939; Fax: 1-310-206-6673.  
\end{document}